\documentclass[review]{elsarticle}

\usepackage{hyperref}
\usepackage{graphicx}
\usepackage{amssymb}
\usepackage{amsmath}
\usepackage{amsthm}

\newtheorem{theorem}{Theorem}
\newtheorem{lemma}{Lemma}
\newtheorem{corollary}{Corollary}

\newtheorem{question}{Question}

\journal{Discrete Applied Mathematics SI: CTW 2020}

\begin{document}

\begin{frontmatter}

\title{Most direct product of graphs are Type~1}

\author[label1]{D. Castonguay}
\address[label1]{INF, Universidade Federal de Goiás, Brazil} 
\author[label2]{C. M. H. de Figueiredo}
\address[label2]{COPPE, Universidade Federal do Rio de Janeiro, Brazil}
\author[label3]{L. A. B. Kowada}
\address[label3]{IC, Universidade Federal Fluminense, Brazil}
\author[label2,label4]{C. S. R. Patr\~ao}
\address[label4]{IFG, Instituto Federal de Goi\'as, Brazil}
\author[label5]{D. Sasaki}
\address[label5]{IME, Universidade do Estado do Rio de Janeiro, Brazil}

\begin{abstract}
A \textit{$k$-total coloring} of a graph $G$ is an assignment of $k$ colors to its elements (vertices and edges) so that adjacent or incident elements have different colors. The total chromatic number is the smallest integer $k$ for which the graph $G$ has a $k$-total coloring. Clearly, this number is at least $\Delta(G)+1$, where $\Delta(G)$ is the maximum degree of $G$. When the lower bound is reached, the graph is said to be Type~1. The upper bound of $\Delta(G)+2$ is a central problem that has been open for fifty years, is verified for graphs with maximum degree 4 but not for regular graphs. 

Most classified direct product of graphs are Type~1.
The particular cases of the direct product of cycle graphs $C_m \times C_n$, 
for $m =3p, 5\ell$ and $8\ell$ 
with $p \geq 2$ and $\ell \geq 1$,
and arbitrary $n \geq 3$, 
were previously known to be Type 1 and 
motivated the conjecture that, 
except for $C_4 \times C_4$, all direct product of cycle graphs $C_m \times C_n$ with $m,n \geq 3$ are Type 1.

We give a general pattern proving that all $C_m \times C_n$ are Type 1, except for $C_4 \times C_4$. 
Additionally, we investigate sufficient conditions  to ensure that the direct product  reaches the lower bound for the total chromatic number. 
\end{abstract}


\end{frontmatter}


\section{Introduction}
\label{intro}

Let $G$ be a simple connected graph with vertex set $V(G)$ and edge set $E(G)$. 
A \textit{$k$-total coloring} of $G$ is an assignment of $k$ colors to its elements (vertices and edges) so that adjacent or incident elements have distinct colors.
The \textit{total chromatic number} $\chi_T(G)$ is the smallest integer $k$ for which $G$ has a $k$-total coloring. Clearly, $\chi_T(G) \geq \Delta(G) +1$, where $\Delta(G)$ is the maximum degree of $G$, and the \textit{Total Coloring Conjecture} (TCC), posed fifty years ago independently by 
Vizing~\cite{vizing} and Behzad et al.~\cite{Behzad},
states that $\chi_T(G) \leq \Delta(G)+2$. 
Graphs with $\chi_T(G)=\Delta(G)+1$ are said to be \textit{Type~$1$} and graphs with $\chi_T(G)=\Delta(G)+2$ are said to be \textit{Type~$2$}. 
In 1977, Kostochka~\cite{kostochka77} verified the TCC for
all graphs with maximum degree 4,
but the TCC has not been settled for all regular graphs, where all vertices have the same degree.
Although the TCC is trivially settled for all bipartite graphs, the problem of determining the total chromatic number of a $k$-regular bipartite graph is NP-hard, for each fixed $k \geq 3$~\cite{ColinSanchez},
exposing how challenging the problem of total coloring is.

The \textit{direct product} (also called \textit{tensor product} or \textit{categorical product}) of two graphs $G$ and $H$ is a graph denoted by $G \times H$, whose vertex set is the Cartesian product of the vertex sets $V(G)$ and $V(H)$
that is $\{(u, v): u \in V(G), v \in V(H)\}$, for which vertices $(u,v)$ and $(u\textsc{\char13},v\textsc{\char13})$
are adjacent if and only if $uu\textsc{\char13} \in E(G)$ and $vv\textsc{\char13} \in E(H)$. The definition clearly implies that  the maximum degree satisfies $\Delta(G \times H)=\Delta(G) \cdot \Delta(H)$, and that the direct product $G \times H$ is a regular graph if and only if both $G$ and $H$ are regular graphs. Concerning the category of graphs, where objects are graphs and morphisms are graph homomorphisms, we know that the direct product $G \times H$ is the categorical product that is defined by projections $p_G: G \times H \rightarrow G$ and $p_H: G \times H \rightarrow H$. The direct product has the commutative property, that is, the graph $G \times H$ is isomorphic to the graph $H \times G$. 
Moreover, 
the direct product $G \times H$ is bipartite if and only if $G$ or $H$ is bipartite, and it is disconnected if and only if $G$ and $H$ are bipartite graphs.
In particular,
in case both $G$ and $H$ are connected bipartite graphs, the direct product $G \times H$ has exactly two bipartite connected components.

A total coloring defines a vertex coloring and an edge coloring, and both coloring problems have been studied with respect to the direct product. A $z$-vertex (resp. edge) coloring  of a graph is an assignment of $z$ colors to its vertices (resp. edges) so that adjacent vertices (resp. incident edges) have distinct colors. The
 chromatic number (resp. index) is the smallest integer $z$ for which a graph has a $z$-vertex (resp. edge) coloring.
Concerning vertex coloring, Hedetniemi conjectured in 1966 that the chromatic number of $G \times H$ would be equal to the minimum of the
chromatic numbers of $G$ and $H$ and recently, fifty years later, the conjecture has been refuted by Shitov~\cite{Shitov}.
Concerning edge coloring, Jaradat~\cite{jaradat} proved that if one factor reaches the lower bound for edge coloring, so does the direct product.

A \textit{cycle graph}, denoted by $C_n, n \geq 3$, is a connected 2-regular graph.
The graph $C_n$ is Type~1 if $n$ is multiple of $3$ and Type~2, otherwise~\cite{yapbook}.
The direct product of cycle graphs $C_m \times C_n$ is a 4-regular graph, 
and it is disconnected precisely when both $m$ and $n$ are even in which case $C_m \times C_n$ consists of two isomorphic 4-regular bipartite connected components each being a spanning subgraph of the complete bipartite graph $K_{\frac{nm}{4},\frac{nm}{4}}$. 
 


 Concerning the total coloring of the direct product, there are few known results.
 Most classified direct product of graphs are Type 1.
 Prnaver and Zmazek~\cite{zmazek} established the TCC for the direct product of a path of length greater or equal to 3 and an arbitrary  graph $G$ with chromatic index $\chi'(G) = \Delta(G)$. They additionally proved, for $m,n \geq 3$, that $\chi_T(P_m \times P_n)$ and 
 $\chi_T(P_m \times C_n)$ are equal to 5.
 Recently, the total chromatic number of direct product of complete graphs has been fully determined as being Type~1 with the exception of $K_2 \times K_2$~\cite{Carol2021}. 
 
  An equitable total coloring is a total coloring where the number of elements colored with each color differs by at most one.
  In 2009, Tong et al.~\cite{tong} showed that the equitable total chromatic number of the Cartesian product between $C_m$ and $C_n$, denoted by $C_{m} \Box C_{n}$, is equal to 5 for all possible values $m, n \geq 3$. It is known that $C_{2n+1} \times C_{2n+1} \cong C_{2n+1} \Box C_{2n+1}$~\cite{geethasum}, therefore we know that $\chi_T(C_{2n+1} \times C_{2n+1})=5$, for all $n\geq 1$.  
 
 In 2018, Geetha and Somasundaram~\cite{geethasum} conjectured that, except for $C_4 \times C_4$, all direct product of cycle graphs $C_m \times C_n$ are Type 1. As evidence, they established three infinite families of Type 1 direct product of cycle graphs:  for arbitrary $n \geq 3$, $\chi_T(C_m \times C_n) = 5$ if $m = 3p, 5\ell, 8\ell$, 
 where $p \geq 2$ and $\ell \geq 1$. 
 In order to describe the claimed total colorings for the three infinite families, they present four tables whose entries are the 5 colors given to suitable matchings between independent sets of vertices that are colored with no conflicts.

In Section~\ref{s:cmcn}, we present a general pattern that gives a 5-total coloring for all graphs $C_m \times C_n$, except for $C_4 \times C_4$. Therefore we ensure that the open remaining infinite families of $C_m \times C_n$ are also Type~1.
In Section~\ref{s:conj}, 
we investigate further conditions that ensure that the direct product $G \times H$ is Type~1. 
We ask whether one factor reaching the lower bound
is enough to ensure that the direct product also reaches the lower bound for the 
total chromatic number.
We manage to classify into Type 1 or Type 2 additional bipartite direct product of graphs.


\section{Total coloring of $C_m \times C_n$}
\label{s:cmcn}

In this section, we prove that the graph $C_m \times C_n$ is Type~1, except for $C_4 \times C_4$. Note that the graph $C_4 \times C_4$ is Type $2$, as it is isomorphic to two copies of $K_{4,4}$, well known to be Type~2, and it is the single exception among the direct product of cycle graphs $C_m \times C_n$.



The present section is devoted to the proof of Theorem~\ref{t: CmCnType1}. 

\begin{theorem}
\label{t: CmCnType1}
Except for $C_4 \times C_4$, the graph $C_m \times C_n$ is Type 1. 
\end{theorem}

\begin{figure}[h]
\centering
\includegraphics[width=0.35\textwidth]{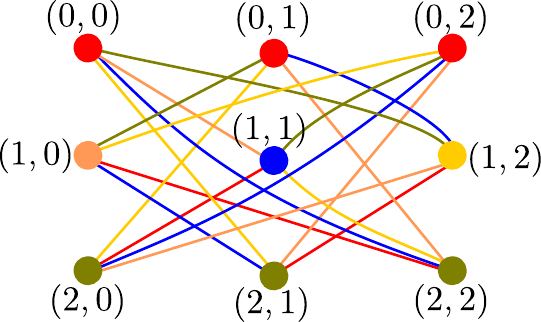}
\caption{A 5-total coloring of $C_3 \times C_3$.}
    \label{C3C3}
\end{figure}

We omit five particular cases that are too small to apply the used technique, but are easy to verify to be Type~1, for instance by using the free open-source mathematics software system Sage Math. They are:
$C_3\times C_3$, $C_3\times C_4$, $C_3\times C_7$, $C_4\times C_7$ and $C_7\times C_7$. Figure \ref{C3C3} presents a 5-total coloring of $C_3 \times C_3$. Therefore, as $C_m \times C_n$ is isomorphic to $C_n \times C_m$, 
we shall consider in our proof $C_m \times C_n$ with $m, n \geq 3$ and $m \neq 3,4,7$. 
We shall write $m = 5k+b$, for $k \geq 0$ and $b = 5,6,8,9$ and $12$. 
Note that as $m \neq 7$, 
the next case is $m = 12$ for which the remainder of the division by 5 is 2.
For instance, to obtain a 5-total coloring of $C_3 \times C_{24}$, we consider the isomorphic graph $C_{24} \times C_3$ and write 
$24 = 5k + b$, with $k = 3$ and $b = 9$.

\begin{figure}[htb] \footnotesize
\centering
\footnotesize
\begin{minipage}[c]{0.44\linewidth}
\centering
\includegraphics[width=0.75\textwidth]{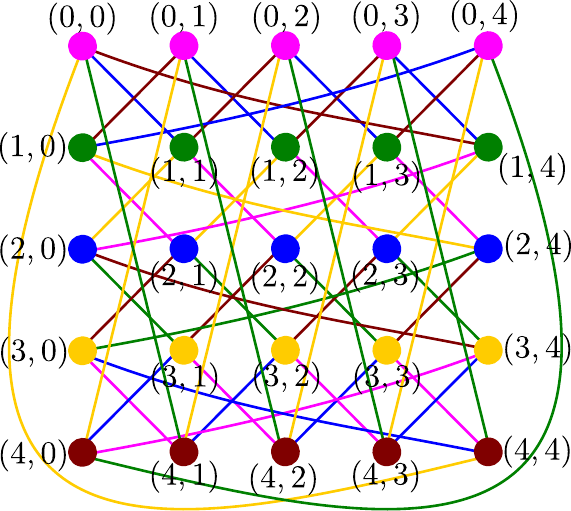}
\end{minipage}
\hfill
\begin{minipage}[c]{0.44\linewidth}
\centering
\includegraphics[width=0.7\textwidth]{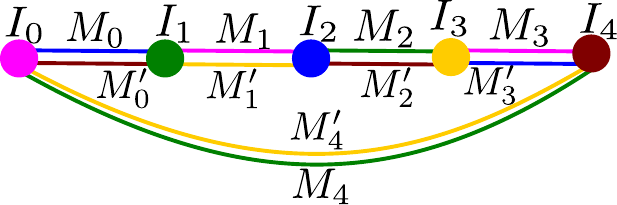}
\end{minipage}
\caption{\footnotesize 
A 5-total coloring of $C_5 \times C_5$ (left) and the respective colored matching quotient $Q[C_5 \times C_5]$ (right).}
\label{C5C5}
\end{figure}

In order to prove Theorem~\ref{t: CmCnType1}, we construct a $5$-total coloring of an auxiliary graph, called matching quotient, from which we obtain a $5$-total coloring of $C_{m} \times C_{n}$.
We use suitable independents sets and matchings between them, inspired by the strategy used by Geetha and Somasundaram~\cite{geethasum} to total color the three particular infinite families.
For $i = 0, \ldots, m-1$, denote by $I_i = \{(i, j)~|~j = 0, \ldots, n-1\}$, $M_i = \{(i, j)((i + 1)~\hbox{mod}~m, (j + 1)~\hbox{mod}~n)~|~j = 0, \ldots,n-1\}$ and $M'_i = \{(i, j)((i + 1)~\hbox{mod}~m, (j-1)~\hbox{mod}~n)~|~j = 0,\ldots, n-1\}$. 
Clearly, each set has $n$ elements, and sets $M_i$ and $M'_i$ are two perfect matchings between independent sets $I_i$ and $I_{i+1}$ in $C_{m} \times C_{n}$.
From that, we define the {\it matching quotient} of $C_m \times C_n$, denoted by $Q[C_m \times C_n]$, as the multigraph where each of its $m$ vertices correspond to an independent set $I_i,  i = 0,\ldots,m-1$, and we have two edges between $I_i$ and $I_{i+1}$ which correspond to $M_i$ and $M_i'$. 
Note that a 5-total coloring of the
matching quotient $Q[C_m \times C_n]$ represents a 5-total coloring of $C_m \times C_n$. 
Figure~\ref{C5C5} presents an example of a matching quotient, by depicting the matching quotient $Q[C_5 \times C_5]$ of $C_5\times C_5$. 

In Subsections \ref{sub:basefamily} and \ref{sub:merge}, we establish a 5-total coloring of the matching quotient of $C_m \times C_n$, proving Theorem \ref{t: CmCnType1}. 
In Subsection \ref{sub:basefamily}, we exhibit a 5-total coloring of the matching quotients of five base infinite families: $C_{5} \times C_{n}, C_{6} \times C_{n}, C_{8} \times C_{n}$, $C_{9} \times C_{n}$ and $C_{12} \times C_{n}$, for $n \geq 3$. 
Note that the base 
infinite families are those where $m =5k+b$ for $k=0$ and $b=5,6,8,9,12$.  We observe that the 5-total
coloring of the base infinite family $C_5 \times C_n$ acts as a pattern.
In Subsection~\ref{sub:merge}, for the matching quotient of $C_{m} \times C_{n}$, with an arbitrary large value of $m$, we observe that we can split this graph into (possibly many) pattern blocks that are identified with the matching quotient of $C_{5} \times C_{n}$, and one base block which is identified with the matching quotient of 
each
base infinite family $C_{5} \times C_{n}$, $C_{6} \times C_{n}$, $C_{8} \times C_{n}$, $C_{9} \times C_{n}$ and $C_{12} \times C_{n}$. 
The 5-total colorings of the matching quotients, given in Subsection \ref{sub:basefamily}, produce a 5-total coloring of each block such that there are no conflicts of colors.
The strategy of splitting the graph into blocks gives
a 5-total coloring of the matching quotient of $C_{m} \times C_{n}$, ensuring that $C_{m} \times C_{n}$ is Type~1.

\subsection{Base infinite families}
\label{sub:basefamily}

\begin{figure}[htb] \footnotesize
\centering
\footnotesize
\begin{minipage}[c]{0.54\linewidth}
\centering
\scalebox{0.8}{
\begin{tabular}{|c|c|c|}
 \hline
   Color & Vertices & Edges \\ \hline
 $1$ (pink) & $I_0$ & $M_1, M_3$\\
\hline
$2$ (green) & $I_1$ & $M_2, M_4$\\
\hline
$3$ (blue) & $I_2$ & $M_0, M'_3$\\
\hline
 $4$ (yellow) & $I_3$ & $M'_1, M'_4$\\
\hline
$5$ (brown) & $I_4$ & $M'_0, M'_2$\\
\hline
\end{tabular}}\end{minipage}
\hfill
\begin{minipage}[c]{0.44\linewidth}
\centering
\includegraphics[width=0.7\textwidth]{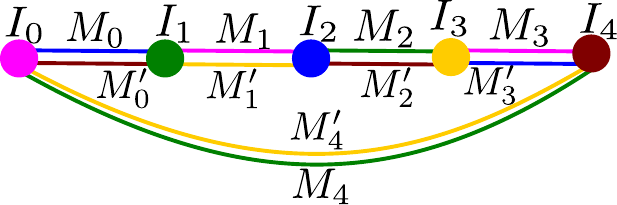}
\end{minipage}
\caption{\footnotesize 
The table and the drawing of a 5-total coloring of the base case $Q[C_5 \times C_n]$.}
\label{t:m=5k}
\end{figure}
\begin{figure}[htb] \footnotesize
\centering
\footnotesize
\begin{minipage}[c]{0.54\linewidth}
\centering
\scalebox{0.8}{
\begin{tabular}{|c|c|c|}
 \hline
   Color & Vertices & Edges \\ \hline
$1$ & $I_0, I_3$ & $M_1, M_4$ \\
 \hline
 $2$ & $I_1, I_4$ & $M_2, M_5$\\
 \hline
 $3$ & $I_2, I_5$ & $M_0, M_3$\\ \hline
 $4$ & $-$ & $M'_
1, M'_3, M'_5$\\
\hline
$5$ & $-$ & $M'_0,M'_2, M'_4$\\
\hline
\end{tabular}}
\end{minipage}
\hfill
\begin{minipage}[c]{0.44\linewidth}
\centering
\includegraphics[width=0.55\textwidth]{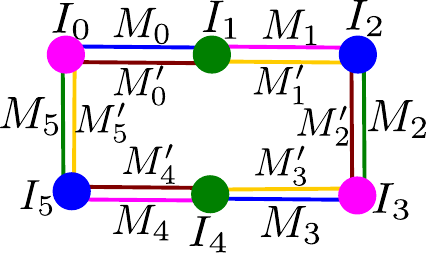}
\end{minipage}
\caption{\footnotesize 
The table and the drawing of a 5-total coloring of the base case $Q[C_6 \times C_n]$.}
\label{t:m=5k+1}
\end{figure}
\begin{figure}[!htb] \footnotesize
\centering
\footnotesize
\begin{minipage}[c]{0.54\linewidth}
\centering
\scalebox{0.89}{
\begin{tabular}{|c|c|c|}
 \hline
   Color & Vertices & Edges \\ \hline
$1$ & $I_0, I_5$ & $M_1, M_3, M_6$\\
\hline
$2$ & $I_1, I_4$ & $M_2, M_5, M_7$\\
\hline
$3$ & $I_2, I_7$ & $M_0, M'_3, M'_5$\\
\hline
$4$ & $I_3, I_6$ & $M'_1, M_4, M'_7$\\
\hline
$5$ & $-$ & $M'_0,M'_2, M'_4, M'_6$\\
\hline
\end{tabular}}
\end{minipage}
\hfill
\begin{minipage}[c]{0.44\linewidth}
\centering
\includegraphics[width=0.65\textwidth]{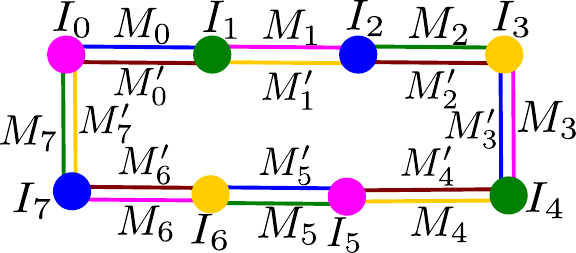}
\end{minipage}
\caption{\footnotesize 
The table and the drawing of a 5-total coloring of the base case $Q[C_8 \times C_n]$.}
\label{t:m=5k+3}
\end{figure}
\begin{figure}[!htb] \footnotesize
\centering
\footnotesize
\begin{minipage}[c]{0.54\linewidth}
\centering
\scalebox{0.85}{
\begin{tabular}{|c|c|c|}
 \hline
   Color & Vertices & Edges \\ \hline
$1$ & $I_0, I_3, I_6$ & $M_1, M_4, M_7$\\
\hline
$2$ & $I_1, I_4, I_7$ & $M_2, M_5, M_8$\\
\hline
$3$ & $I_2$ & $M_0, M_3, M'_5, M'_7$\\
\hline
$4$ & $I_5$ & $M'_1, M'_3, M_6, M'_8$\\
\hline
$5$ & $I_8$ & $M'_0, M'_2, M'_4, M'_6$\\
\hline
\end{tabular}}
\end{minipage}
\hfill
\begin{minipage}[c]{0.44\linewidth}
\centering
\includegraphics[width=0.7\textwidth]{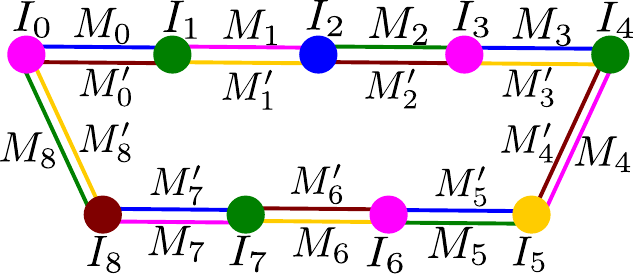}
\end{minipage}
\caption{\footnotesize 
The table and the drawing of a 5-total coloring of the base case $Q[C_9 \times C_n]$.}
\label{t:m=5k+4}
\end{figure}
\begin{figure}[!htb] \footnotesize
\centering
\footnotesize
\begin{minipage}[c]{0.54\linewidth}
\centering
\scalebox{0.85}{
\begin{tabular}{|c|c|c|}
 \hline
   Color  & Vertices & Edges \\ \hline
$1$ & $I_0, I_3, I_6, I_9$ & $M_1, M_4, M_7, M_{10}$\\
\hline
$2$ & $I_1, I_8$ & $M_2, M'_4, M_6, M_9, M_{11}$\\
\hline
$3$ & $I_2, I_7$ & $M_0, M_3, M_5, M_8, M'_{10}$\\
\hline
$4$ & $I_5, I_{10}$ & $M'_1, M'_3, M'_6, M'_8, M'_{11}$\\
\hline
$5$ & $I_4, I_{11}$ & $M'_0, M'_2, M'_5, M'_7, M'_9$\\
\hline
\end{tabular}}
\end{minipage}
\hfill
\begin{minipage}[c]{0.44\linewidth}
\centering
\includegraphics[width=0.75\textwidth]{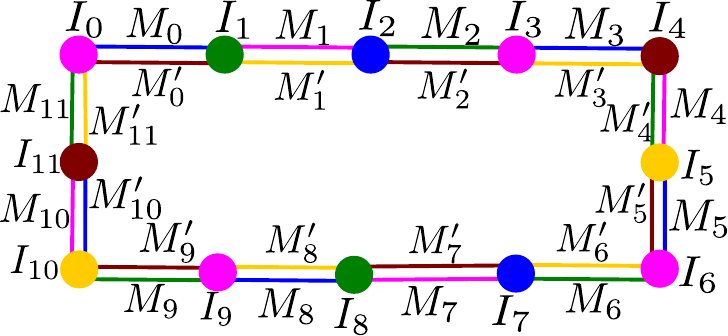}
\end{minipage}
\caption{\footnotesize 
The table and the drawing of a 5-total coloring of the base case $Q[C_{12} \times C_n]$.}
\label{t:m=5k+2}
\end{figure}

We consider first the base infinite families $Q[C_{m} \times C_n]$ with $m=5,6,8,9,12$ and $n \geq 3$. 
We refer to Figures~\ref{t:m=5k},~\ref{t:m=5k+1},~\ref{t:m=5k+3},~\ref{t:m=5k+4} and~\ref{t:m=5k+2} for the 5-total coloring of each base case. Note that, the 5-total colorings of the base infinite families have important features in common: the same color~1 (pink) given to the vertex $I_0$, the same color~2 (green) given to the matching $M_{m-1}$ and the same color~4 (yellow) given to the matching $M'_{m-1}$. 
These shared features provide the needed compatibility that allows us to define a common pattern used when we deal with larger values of $m$.


\subsection{Merging the pattern to generate a 5-total coloring for arbitrary $m$}
\label{sub:merge}

To obtain a 5-total coloring of the matching quotient $Q[C_m \times C_n]$ for an arbitrary large value of $m$, 
we repeatedly merge the pattern block given by 5-total coloring of the matching quotient $Q[C_5 \times C_n]$ with the 5-total coloring of its base block $Q[C_b \times C_n]$, for $b=5,6,8,9,12$. Note that the colored $Q[C_5 \times C_n]$ is the first base case and is also the only pattern used for an arbitrary value of $m$ independently of its base case. 

Recall that, as we argued in the beginning of Section~\ref{s:cmcn}, by swapping $m$ and $n$, we are always able to consider $n \geq 3$ and write a large value of $m \geq 10$ as
$m=5k+b$, for $k \geq 1$ and $b=5,6,8,9,12$.
We optimally color first its base block $C_b \times C_n$ and then repeatedly merge with $k$ copies of the optimally colored pattern block $C_5 \times C_n$.
So the 5-total coloring of $Q[C_m \times C_n]$ is defined by two steps as follows:


\begin{itemize}

\item \textbf{Base step}: For each $i=0,  \ldots, b-1$, the color of $I_i$ (respectively, $M_i$ and $M'_i$) in $Q[C_m \times C_n]$ is the same as the color of $I_i$ (respectively, $M_i$ and $M'_i$) in its base case $Q[C_b \times C_n]$.%



\item  \textbf{Pattern step}: For each $i= b,  \ldots m-1$, write $t=(i-b)~\hbox{mod}~5$, and the color of $I_i$ (respectively, $M_i$ and $M'_i$) in $Q[C_m \times C_n]$ is the same as the color of $I_{t}$ (respectively, $M_{t}$ and $M'_{t}$) in the pattern $Q[C_5\times C_n]$.
\end{itemize}

\begin{figure}[!htb]
    \centering
    \includegraphics[width=0.5\textwidth]{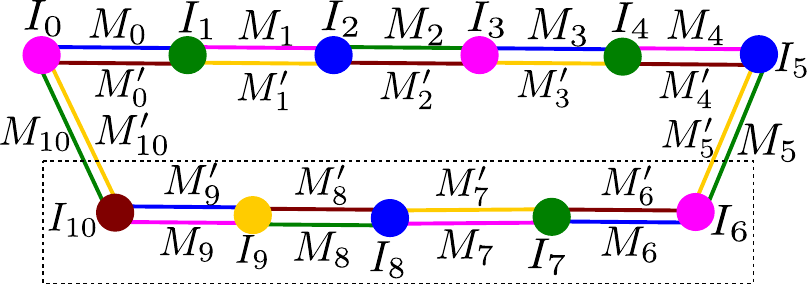}
    \caption{\footnotesize A 5-total coloring of the matching quotient $Q[C_{11} \times C_n]$. This 5-total coloring is obtained by merging once the 5-total coloring of the highlighted pattern $Q[C_{5} \times C_n]$ into the 5-total coloring of $Q[C_{6} \times C_n]$.}
    \label{fig:C11Cn}
\end{figure}

\begin{figure}[htb]
    \centering
    \includegraphics[width=0.8\textwidth]{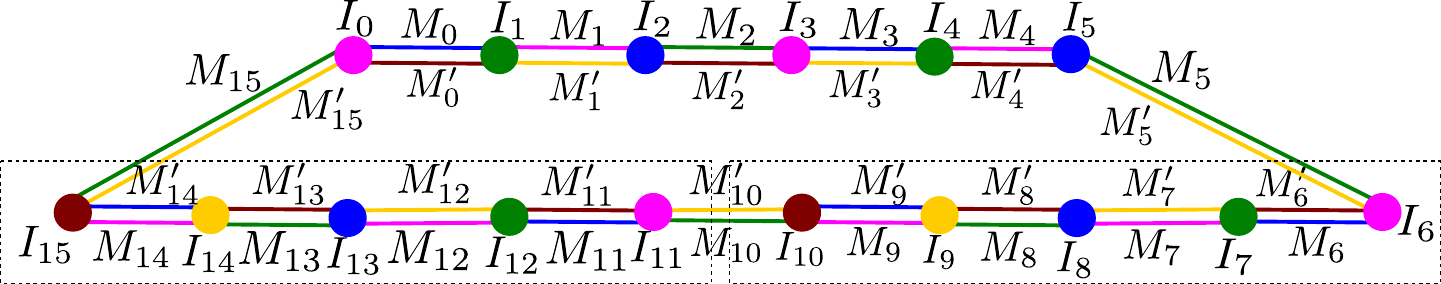}
    \caption{\footnotesize A 5-total coloring of the matching quotient $Q[C_{16} \times C_n]$. This 5-total coloring is obtained by merging twice the 5-total coloring of the highlighted pattern $Q[C_{5} \times C_n]$ into the 5-total coloring of $Q[C_{6} \times C_n]$.}
    \label{fig:C16Cn}
\end{figure}

For instance, consider $m=11$ and please refer to Figure ~\ref{fig:C11Cn}.
Note that, in the base step, we color the elements $I_i, M_i$ and $M'_i$, for $i=0, \ldots, 5$, of $Q[C_{11} \times C_n]$ with the same colors as its base infinite family $Q[C_6 \times C_n]$. Now, in the pattern step, we color the elements $I_i, M_i$ and $M'_i$, for $i=6,7,8,9,10$ of $Q[C_{11} \times C_n]$ with the same colors as the pattern $Q[C_5 \times C_n]$ (as in Figure~\ref{t:m=5k}).
Analogously, when $m=16$ we merge the pattern $Q[C_5 \times C_n]$ twice into $Q[C_6 \times C_n]$ to obtain a 5-total coloring of the matching quotient $Q[C_{16} \times C_n]$ as highlighted in Figure~\ref{fig:C16Cn} by elements $I_i, M_i$ and $M'_i$, for $i=6, \ldots, 15$. Thus, for a general $m= 5k+b$ we merge $k$ patterns $Q[C_5 \times C_n]$ into the corresponding base infinite family $Q[C_b \times C_n]$ to obtain a 5-total coloring of the matching quotient $Q[C_{5k+b} \times C_n]$. 



Note that there  is  no  conflict  between  the  assigned  colors  in  the defined 5-total coloring for an arbitrary value of $m$. 
Indeed, we already know that each base infinite family has color 1 (pink) to $I_0$, color 2 (green) to $M_{b-1}$ and color 4 (yellow) to $M'_{b-1}$. Note that, regardless of how many times we use the pattern, there is no conflict between the patterns, as the edges colored 2 (green) and 4 (yellow) in the case base $Q[C_5 \times C_n]$ are the ones used between the patterns. Also, there is no conflict between the base case and the pattern, as the edges colored 2 (green) and 4 (yellow) between $I_{b-1}$ and $I_{b}$ (pink), and between $I_{0}$ colored 1 (pink) and $I_{m-1}$, are both used in our base cases, where both $I_{b-1}$ and $I_{m-1}$ have a colored 1 (pink) neighbor vertex.



 
 



\vspace{0.2cm}



\section{On total coloring bipartite direct product of graphs}
\label{s:conj}

In this section, we propose and investigate two questions motivated by the search for a general pattern for the classification into Type~1 or Type~2 of the direct product of two graphs. In this sense, it is natural to seek for sufficient conditions for the direct product to be Type 1. Prnaver and Zmazek~\cite{zmazek} previously proved that if $G$ admits a $\Delta(G)$-edge coloring, then $G\times P_m$, for $m \geq 3$, is Type~1. Mackeigan and Janssen~\cite{JaMa19} subsequently proved that if $G \times K_2$ is Type 1, then $G \times H$ is also Type 1, for any bipartite graph $H$. Recall that for edge coloring, Jaradat~\cite{jaradat} proved that if one factor reaches the lower bound for edge coloring, so does the direct product. We investigate whether an analogous property holds for total coloring:



\begin{question}
\label{quest-type2}
Given a Type~1 graph $G$ and an arbitrary graph $H$, is the direct product $G \times H$ Type~1?
\end{question}

The analogous question has been considered for the Cartesian product, but has only been partially answered. It is known that if the factor with largest vertex degree is of Type 1, then the Cartesian product is also of Type 1~\cite{zmazek2002}.
 So far, all known Type~2 direct product of two graphs are the direct product $G \times H$, where $G$ and $H$ are Type~2, including cases with  $G = H$. 
 The known Type~2 direct product of two graphs are:
 $K_2 \times K_2$, $C_4 \times C_4$, $K_{n,n}\times K_{m,m}$,
 and $C_m \times K_2$ for $m$ not a multiple of 3.
 On the other hand, a Type~1 direct product of two graphs $G \times H$ can be obtained when $G$ and $H$ are Type~1, when $G$ and $H$ are Type~2, or else when one of them is Type~1 and the other is Type~2. 
 For instance, $K_m$ is Type~1 when $m$ is odd and Type~2 when $m$ is even, and yet the direct product $K_m \times K_n$ is Type~1 when both $m$ and $n$ are odd~\cite{Carol2021}, when $m,n\neq2$ are both even~\cite{geethasum}, or else when $m$ or $n$ is even~\cite{JaMa19}; whereas when $m=n=2$, the graph $K_2 \times K_2$ is Type~2. 
 The present work established that the direct product $C_m \times C_n$ is Type~1 when $m,n\neq4$, and yet $C_m$ is Type~1 when $m$ is multiple of 3 and Type~2 otherwise; whereas when $m=n=4$, the graph $C_4\times C_4$ is Type~2.  



We contribute to Question \ref{quest-type2} by giving positive evidences. A regular graph $G$ is \textit{conformable} if $G$ admits a vertex coloring with $\Delta(G) + 1$ colors such that the number of vertices in each color class has the same parity as $|V(G)|$~\cite{chetwyndhiltoncheng}. It is known that every Type 1 graph must satisfy the conformable condition. The converse is not true, but being conformable helps to identify whether a graph has the potential to be Type 1 or to be sure that it cannot be Type 1.
In Theorem \ref{t:conformable}, we show a sufficient condition on the graph $G$ for the direct product of regular graphs $G\times H$ to be conformable.

\begin{theorem}
\label{t:conformable}
Let $G$ and $H$ be two regular graphs. If $G$ is conformable, then $G \times H$ is conformable. 
\end{theorem}

\begin{figure}[!ht]
\centering
\includegraphics[scale=0.65]{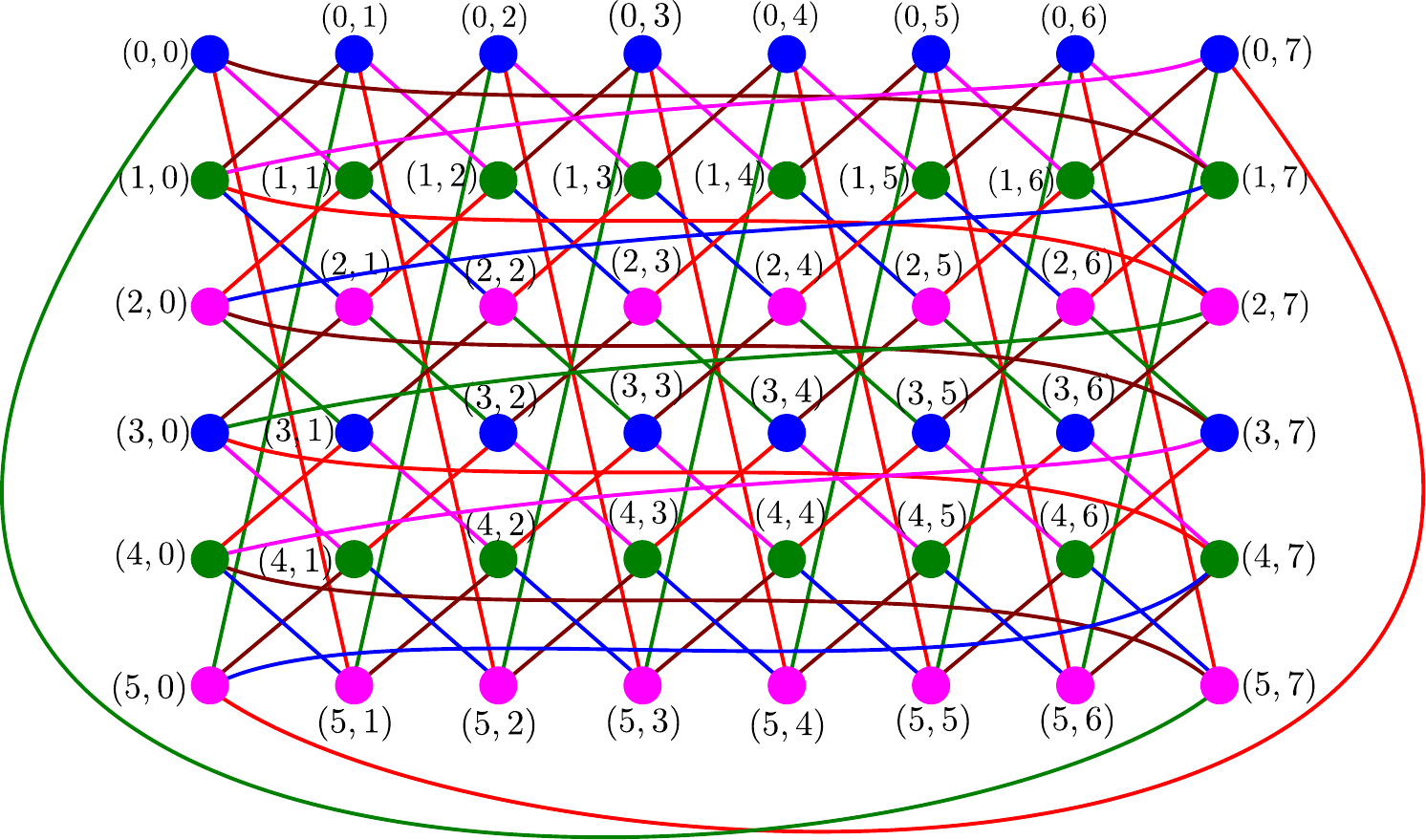}
\caption{A 5-total coloring of $C_6 \times C_8$ which is also the conformable coloring of Theorem~\ref{t:conformable}.}
\label{C6C8}
\end{figure}

\begin{proof}
Since $G$ is by hypothesis conformable, let us consider a vertex coloring $f: V(G) \to \{1, \ldots, \Delta(G)+1\}$ such that, for every $i=1, \ldots, \Delta(G)+1$, the color class $F_i=f^{-1}(i)$ has cardinality of the same parity as $\vert V(G) \vert$. 
Consider one of the projections that define the direct product 
$p:V(G \times H) \to V(G)$. Therefore, we have a function $f \circ p :V(G \times H) \to \{1, \ldots, \Delta(G)+1\}$, which is a vertex coloring of $G \times H$ such that every color class consists of the vertices in the Cartesian product 
of the sets of 
vertices of $F_i$ and $V(H)$, denoted by $F_i \times V(H)$, for $i=1, \ldots, \Delta(G)+1$. 

We consider two cases. First, consider the case when $G \times H$ is a graph of even order. In this case, $\vert V(G) \vert$ or $\vert V(H) \vert$ is even. Recall that $\Delta(G \times H) = \Delta(G) \cdot \Delta(H)$. Consider a vertex coloring $\gamma: V(G \times H) \to \{1, \ldots, \Delta(G) \cdot \Delta(H)+1\}$, defined by $\gamma(x,a) = f(x)$. Note that this function is actually obtained from $f \circ p$ by changing the codomain, thus this is also a vertex coloring of $G \times H$, but possibly there are empty color classes. An empty set has cardinality zero, which is of even parity. We have to prove that $\gamma$ is conformable, that is, we have to prove that every non-empty color class of $\gamma$ has an even cardinality. But, a non-empty color class of $\gamma$ is of the form $F_i \times V(H)$,
where $F_i$ a color class of $f$. Since $F_i$ has the same cardinality of $\vert V(G) \vert$ and $\vert F_i \times V(H)\vert = \vert F_i \vert \cdot \vert V(H) \vert$, if $\vert V(G) \vert$ or $\vert V(H) \vert$ is even, then $\vert F_i \times V(H)\vert$ is also even, and thus, $\gamma$ is a conformable coloring of $G \times H$.

Now, consider the case when $G \times H$ is a graph of odd order, that is $\vert V(G) \vert$ and $\vert V(H) \vert$ are odd. Since $G$ and $H$ are regular graphs, we have that $\Delta(G)$ and $\Delta(H)$ are even. 
Observe that all color classes of $c \circ p$ have odd cardinality. Next, we define a conformable coloring $\gamma: V(G \times H) \rightarrow \{1, \ldots, \Delta(G) \cdot \Delta(H)+1\}$.  Observe that all color classes of $\gamma$ have odd cardinality, hence they must have at least one vertex. The idea is to remove even amount of vertices of each color class of $f \circ p$. In this way, the parity of these color classes is preserved and additionally we define new color classes of cardinality one, for each of this removed vertices.
Consider $V(H) = \{v_1, \ldots, v_{n}\}$ and consider, for each $i=1, \ldots, \Delta(G)+1$, a fixed element $u_i$ of $F_i$. 
We define:
$$A_i = \left\{\begin{array}{ll} \{(u_i, v_j) \mid j=1, \ldots, \Delta(H)-2 \}, \text{if $i=1, \ldots, \Delta(G)/2$} \\ \{(u_i, v_j) \mid j=1, \ldots, \Delta(H) \}, \text{if $i=(\Delta(G)/2)+1, \ldots, \Delta(G)$} \end{array}\right.$$

 In addition, define $A=\bigcup\limits_{i=1}^{\Delta(G)} A_i$.  We construct a vertex coloring $\gamma$ of $G \times H$  by defining its color classes, each one is a subset of a color class of $f \circ p$. Each color class of $\gamma$ is $F_i \times V(H) - A_i$, for $i=1, \ldots, \Delta(G)$, $F_{\Delta(G)+1}$ and $\{(u,v)\}$, for $(u,v) \in A$. In order to show that $\gamma$ is conformable, we have to show that $\gamma$ has $\Delta(G) \cdot \Delta(H)+1$ color classes and that each color class has an odd cardinality, as follows. 

First, the number of color classes of $\gamma$ is $\Delta(G) + 1 + \vert A \vert$.  Since 

 
$$\vert A \vert = \sum\limits_{i=1}^{\Delta(G)} \vert A_i \vert
= (\sum\limits_{i=1}^{\frac{\Delta(G)}{2}} \vert A_i \vert) + (\sum\limits_{i=(\frac{\Delta(G)}{2})+1}^{\Delta(G)} \vert A_i \vert)= \frac{\Delta(G)}{2}(\Delta(H)-2) + \frac{\Delta(G)}{2}\Delta(H)$$ $$=(\frac{\Delta(G)}{2})(2\Delta(H)-2)
= \Delta(G)(\Delta(H)-1),$$


we have that $\Delta(G) + 1 + \vert A \vert$ = $\Delta(G) \cdot \Delta(H)+1$.

Finally, we prove that each color class of $\gamma$ has an odd cardinality. Recall that $H$ is of odd order and has even degree. Since $f$ is a conformable coloring of $G$ and $G$ is of odd order, each $F_i$ is of odd cardinality, for  $i=1, \ldots, \Delta(G)+1$. For each $i=1, \ldots, \Delta(G)/2$, $F_i \times V(H) - A_i$ has odd cardinality $\vert F_i \vert \cdot \vert V(H) \vert - \Delta(H)+2$. Similarly, for each $i=(\Delta(G)/2)+1, \ldots, \Delta(G)$, $F_i \times V(H) - A_i$ has odd cardinality $\vert F_i \vert \cdot \vert V(H) \vert - \Delta(H)$. Clearly, $F_{\Delta(G)+1}$ and $\{(u,v)\}$, for $(u,v) \in A$ have odd cardinality. Therefore, $\gamma$ is a conformable coloring of  $G \times H$.

For examples of the odd and even cases, see Figures~\ref{C3C3} and~\ref{C6C8}, respectively. Note that the conformable colorings of the base infinite families in Section \ref{s:cmcn} are not obtained in the same way, except $C_6\times C_n$, $n\geq 3$. 
\end{proof}

By Theorem \ref{t:conformable}, we contribute to Question \ref{quest-type2} since given regular graphs $G$ and $H$ with $G$ of Type 1, we know that $G\times H$ is conformable. 
Conformable graphs of odd order and sufficient large maximum degree are Type~1, see Chew~\cite{chew}. For example, in \cite{Carol2021}, we use this fact together with Hamiltonian decompositions, to give a full classification of the total chromatic number of the direct product of complete graphs $K_m \times K_n$.



Next 
lemma presents a sufficient condition on the graph $G$ for the direct product $G \times K_2$ to be Type 1, which leads to a corollary that answers Question~1 positively when one factor is Type 1 and the other is bipartite.




\begin{lemma}
\label{t: GxK2}
If $G$ is Type~1, then $G\times K_2$ is Type 1.
\end{lemma}

\begin{proof}
Let $f:V(G) \cup E(G) \rightarrow \{1, \cdots, \Delta(G) + 1\}$ be a total coloring of a Type~1 graph $G$. 
Consider one of the projections that define the direct product 
$p: V(G \times K_2) \cup E(G \times K_2) \rightarrow V(G) \cup E(G)$ and the composite function $f \circ p: V(G \times K_2) \cup E(G \times K_2) \rightarrow \{1, \cdots , \Delta (G)+1\}$. Observe that $\Delta (G \times K_2)=\Delta (G)$. 



First note that, as before, $f \circ p$, when restricted to the vertices of $G \times K_2$, is a vertex coloring.
Second, let $(x,i)(y,j) \in E(G \times K_2)$ and suppose this edge $(x,i)(y,j)$ and its endvertex $(x,i)$ have the same color, that is $(f \circ p)((x,i)(y,j)) = (f \circ p)(x,i)$. Thus, $f(xy) = f(x)$, a contradiction since $x$ is an endvertex of the edge $xy$ in $G$.
Finally, let $(x,i)(y,j)$ and $(y,j)(z,k)$ be two adjacent edges of $G \times K_2$ and suppose that these edges have the same colors, that is $(f \circ p)((x,i)(y,j))= (f \circ p)((y,j)(z,k))$. Thus, $f(xy) = f(yz)$, a contradiction since $xy$ and $yz$ are adjacent edges in $G$.
\end{proof}

\begin{figure}[!ht]
\centering
\includegraphics[scale=0.7]{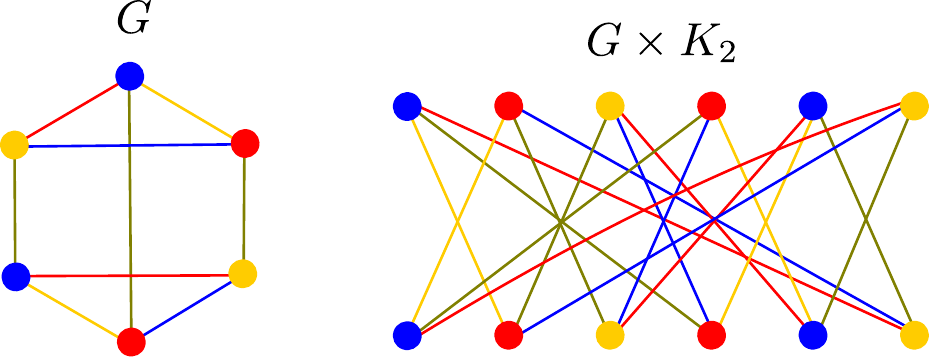}
\caption{A 4-total coloring of $G \times K_2$ obtained by a  4-total coloring of $G$.}
\label{GK2}
\end{figure}

Figure \ref{GK2} presents a 4-total coloring of $G \times K_2$, where $G$ is the Cartesian product $G=C_3\Box K_2$.
The above 
lemma along with the above mentioned result of Mackeigan and Janssen~\cite{JaMa19} give the following corollary:


\begin{corollary}
\label{cor}

If $G$ is Type~1 and $H$ is bipartite, then $G \times H$ is Type~1. 
\end{corollary}

We remark that $C_3 \times C_4$ is Type~1, which agrees with Corollary~\ref{cor}, as $C_3$ is Type~1 and $C_4$ is a bipartite graph.
We remark that $C_m \times K_2$ is Type~1 if and only if $m$ is multiple of~3.
The converse of Corollary~\ref{cor} is not true, since there are many examples of a Type~2 graph $G$ such that the direct product of $G \times K_2$ is Type~$1$. For instance, $K_m$ for even $m$ is Type~2, and yet $K_m \times K_2$ is the complete bipartite graph $K_{m,m}$ minus a perfect matching, known to be Type~1 for $m\geq4$~\cite{yapbook}. 

Also, for two complete bipartite graphs $K_{m,m'}$ and $K_{n,n'}$,
the direct product $K_{m,m'} \times K_{n,n'}$ is Type~2 if and only if $m=m'$ and $n=n'$.
Otherwise, it is Type~1. Indeed, note that if $m\neq m'$ or $n\neq n'$, then $K_{m,m'} \times K_{n,n'}$ is Type~1, by Corollary~\ref{cor} since it is known that $K_{m,m'}$ is Type~1~\cite{yapbook}.
On the other hand, if $m=m'$ and $n=n'$, the graph $K_{m,m} \times K_{n,n}$ is isomorphic to two copies of $K_{mn,mn}$~\cite{jha}. It is known that $K_{mn,mn}$ is Type~2~\cite{yapbook} and thus also $K_{m,m} \times K_{n,n}$.

If $G$ is bipartite and Type~1, then Corollary~\ref{cor} implies that $G \times G$ is Type~1 as well.
In this context, we conclude by proposing the property for a general Type~1 graph:


\begin{question}
\label{quest}
Given a Type~1 graph $G$, is the direct product $G \times G$ Type~1 as well?
\end{question}

\textbf{Acknowledgements}

This work is partially supported by the Brazilian agencies CNPq (Grant numbers: 
302823/2016-6, 407635/2018-1 and 313797/2020-0) and FAPERJ (Grant numbers: CNE E-26/202.793/2017 and ARC E-26/010.002674/2019).




\end{document}